\documentclass[a4paper,11pt,twoside]{article}

\usepackage{geometry}
\usepackage[utf8]{inputenc}
\usepackage{lmodern}
\usepackage[T1]{fontenc}
\usepackage{amsmath,amssymb,amsthm,empheq,cases}
\usepackage{hyperref}
\hypersetup{colorlinks,
            citecolor=red, 
            filecolor=black,
            linkcolor=blue,
            urlcolor=black}
\usepackage{tikz}
\usepackage{graphicx}
\usepackage{url}

\def \dis {\displaystyle}

\def \NN {\mathbb N}

\def \RR {\mathbb R}

\def \A {\mathcal{A}}
\def \B {\mathcal{B}}

\def \F {\mathcal{F}}

\def \L {\mathcal{L}}

\def \P {\mathcal{P}}

\def \ecart {\noalign{\medskip}}

\theoremstyle{definition}
\newtheorem{Th}{Theorem}[section]

\newtheorem{Rem}[Th]{Remark}
 
\def \refs #1{Section~\ref{#1}}

\def \refT #1{Theorem~\ref{#1}}

\def \refR #1{Remark~\ref{#1}}

\title{Generalized diffusion problems in a conical domain, part I}
\author{Rabah Labbas, St\'ephane Maingot \& Alexandre Thorel \\ \ecart
\scriptsize R. L., S. M. \& A. T.: Normandie Univ, UNIHAVRE, LMAH, FR-CNRS-3335, ISCN, 76600 Le Havre, France. \\ \ecart 
\scriptsize rabah.labbas@univ-lehavre.fr, stephane.maingot@univ-lehavre.fr, alexandre.thorel@univ-lehavre.fr}
\date{}

\begin{document}

\maketitle

\begin{abstract}
The purpose of this article (composed of two parts) is the study of the generalized dispersal operator of a reaction-diffusion equation in $L^p$-spaces set in the finite conical domain $S_{\omega,\rho}$ of angle $\omega>0$ and radius $\rho>0$ in $\RR^2$. 

 This first part is devoted to the behaviour of the solution near the top of the cone which is completely described in the weighted Sobolev space $W^{4,p}_{3-\frac{1}{p}}(S_{\omega,\rho})$, see \refT{Th principal}. 
\end{abstract}

\section{Introduction}

This work required the use of many calculations and non-trivial checks linked, among other, to the theory of sums of linear operators. This is why we were forced to split this work into two more or less independent parts.

In this first part, we consider the generalized following reaction-diffusion equation
\begin{equation}\label{eq evolution}
\left\{\begin{array}{l}
\dis\frac{\partial u}{\partial t} = -\Delta^2 u + k \Delta u + f^*(u) \quad \text{in} ~ \RR_+\times\Omega, \\ \ecart
u(0) = u_0 \text{ given} \\
\text{Boundary conditions for } u \text{ on } \partial\Omega,
\end{array}\right.
\end{equation}
where $k$ is a positive number, $f^*$ is a non-linear reaction function and $\Omega$ is a bounded conical domain. This work is a natural continuation of the one studied in \cite{LMMT}. The originality of this work lies in the fact that the open set $\Omega$ is conical whereas in \cite{LMMT}, it was cylindrical.

The study of the spatial operator, that is the linear combination of the laplacian and the bilaplacian, in \eqref{eq evolution} required the analysis of a sum of unbounded linear operators in a Banach space, carried out in part II, see \cite{Cone P2}. Note that, by similar techniques, we can prove that the dispersal operator in \eqref{eq evolution} generates an analytic semigroup, see for instance \cite{LMT}.

Such problems, set in conical domains, model many concrete situations related to pollution for instance.

More precisely, in this work, we consider the following 2 dimensional sector
\begin{equation*}
\Omega = S_{\omega }=\left\{ (x,y)=(r\cos \theta ,r\sin \theta ):r>0\text{ and }%
0<\theta <\omega \right\} ,
\end{equation*}
with its lateral edges 
\begin{equation*}
\left\{ 
\begin{array}{l}
\Gamma _{0}= (0,+\infty) \times \left\{ 0\right\}  \\ 
\Gamma _{\omega }=\left\{ (r\cos \omega ,r\sin \omega ):\text{ }r>0\text{ }\right\} ;
\end{array}\right. 
\end{equation*} 
here $\omega \in (0, 2\pi]$.

It is known that the study of problem \eqref{eq evolution} in $L^p$-spaces needs the analyze of the following linear stationary problem 
\begin{equation}\label{Pb cone infini}
\left\{ \begin{array}{ll}
\Delta ^{2}u-k\Delta u=f &\text{in }S_{\omega } \\ \ecart 
u=\dfrac{\partial u}{\partial n}=0 & \text{on }\Gamma _{0}\cup \Gamma
_{\omega }.
\end{array}\right. 
\end{equation}

This work is inspired by the one done in \cite{geymonat-grisvard} and \cite{grisvard}, where the author has considered, in a hilbertian framework, the following boundary problem
\begin{equation*}\label{Pb Grisvard}
\left\{ \begin{array}{ll}
\Delta^2 u = 0 &\text{in }S_{\omega } \\ \ecart 
u=\dfrac{\partial u}{\partial n}=0 & \text{on }\Gamma _{0}\cup \Gamma
_{\omega }.
\end{array}\right. 
\end{equation*}
The author has proved that the solution of this problem writes as a "superposition" of particular solutions with separate variables of the form $\chi_{1,j}(r)\chi_{2,j}(\theta)$, for $j \in \NN$. The basic tools used are based on the compact operators belonging to the so-called Carleman class and the Fredholm determinants. 

Let us recall some known results concerning the biharmonic equation in a conical domain or in a Lipschitz domain. In \cite{pipher}, the authors gave many estimates concerning the solution of the Dirichlet problem in $L^p$ for the biharmonic equation in Lipschitz domain. In \cite{barton}, many results are given for general higher-order elliptic equations in non smooth-domains. In \cite{tami}, the author has studied the following problem
\begin{equation*}
\left\{ \begin{array}{ll}
\Delta^2 u = f &\text{in }S_{\omega,1} \\ \ecart 
u=\Delta u =0 & \text{on } \partial S_{\omega,1},
\end{array}\right. 
\end{equation*}
where $ f \in L^2(S_{\omega,1})$. He has proved the two following results
\begin{enumerate}
\item If $\omega < \pi$, the variational solution writes, in the neighbourhood of $O$, as 
$$u_\omega = u_{1,\omega} + u_{2,\omega} + u_{3,\omega},$$
with $u_{1,\omega} \in H^{1+\frac{\pi}{\omega}-\varepsilon}$, $u_{2,\omega} \in H^{2+\frac{\pi}{\omega}-\varepsilon}$ and $u_{3,\omega} \in H^{4}$, for a small $\varepsilon > 0$.  

\item If $\omega = \pi$, in the neighbourhood of $O$, the variational solution is verifies  
$$u_\pi \in H^4.$$
\end{enumerate}

This article is organized as follows. In \refs{Sect Statement of result}, we state our problem in polar coordinates and our main result in \refT{Th principal}. In \refs{Sect Reformulation} we reformulate our problem as a sum of linear operators. Then, \refs{Sect proof of main Th} is devoted to its proof.

\section{Statement of the main result}\label{Sect Statement of result}

We introduce the following polar variables function  
\begin{equation*}
v(r,\theta )=u(r\cos \theta ,r\sin \theta )=u(x,y).
\end{equation*}
It is known that the laplacian and the bilaplacian, in polar coordinates, respectively write
\begin{equation}\label{Lambda 1}
\Delta u = \frac{1}{r^{2}}\left[ \left( r\dfrac{\partial }{\partial r}\right)
^{2}+\dfrac{\partial ^{2}}{\partial \theta ^{2}}\right] v=\dfrac{\partial
^{2}v}{\partial r^{2}}+\frac{1}{r}\dfrac{\partial v}{\partial r}+\frac{1}{%
r^{2}}\dfrac{\partial ^{2}v}{\partial \theta ^{2}} := \Lambda_1 v,
\end{equation}
and 
\begin{equation}\label{Lambda 2}
\begin{array}{rcl}
\Delta ^{2}u &\hspace{-0.1cm}=&\hspace{-0.1cm} \dis\left( \dfrac{\partial ^{2}}{\partial r^{2}}+\frac{1}{r}\dfrac{\partial }{\partial r}+\frac{1}{r^{2}}\dfrac{\partial ^{2}}{\partial \theta^{2}}\right) ^{2}v := \Lambda_2 v \\ \ecart
&=&\hspace{-0.1cm}\dis \dfrac{\partial ^{4}v}{\partial r^{4}}+\frac{2}{r^{2}}\dfrac{\partial
^{4}v}{\partial r^{2}\partial \theta ^{2}} +\frac{1}{r^{4}}\dfrac{\partial ^{4}v}{\partial \theta ^{4}} + \frac{2}{r}\dfrac{\partial^{3}v}{\partial r^{3}} -\frac{2}{r^{3}}\dfrac{\partial ^{3}v}{\partial r\partial \theta ^{2}}-\frac{1}{r^{2}}\dfrac{\partial ^{2}v}{\partial r^{2}}+\frac{4}{r^{4}}\dfrac{\partial ^{2}v}{\partial \theta ^{2}}+\frac{1}{r^{3}}\dfrac{\partial v}{\partial r} .
\end{array}
\end{equation}

\begin{Rem}
We can generalize this work to the dimension $n$ :
\begin{equation*}
\Delta u=\dfrac{\partial ^{2}v}{\partial r^{2}}+\frac{n-1}{r}\dfrac{\partial
v}{\partial r}+\frac{1}{r^{2}}\Delta ^{\prime }v,
\end{equation*}
where $\Delta ^{\prime }$ is the Laplace-Beltrami operator.
\end{Rem}
We set
$$f(x,y)=f\left( r\cos \theta ,r\sin \theta
\right) =g(r,\theta ).$$
Let $\gamma \in \RR$. In the sequel, we will use the following weighted space Sobolev spaces on $S_{\omega}$ by:
\begin{equation}\label{Def Wmp}
W^{m,p}_\gamma (S_{\omega}) = \left\{v \in \mathcal{D}'(S_{\omega}), \forall\, (i,j) \in \mathbb{N}^2 : 0 \leqslant i+j \leqslant m, ~ r^{\gamma-m+i+j} \frac{\partial^{i+j} v}{\partial r^i \partial \theta^j} \in L^p(S_{\omega}) \right\}.
\end{equation}
Note that $f \in L^p(\Omega)$ means that $g \in L^p_{\frac{1}{p}} (S_\omega) := W^{0,p}_{\frac{1}{p}}(S_\omega)$. In fact, we have
$$\int_\Omega |f(x,y)|^p dx dy = \int_{S_\omega} |g(r,\theta)|^p r dr d\theta = \int_{S_\omega} |r^{\frac{1}{p}}g(r,\theta)|^p dr d\theta.$$
We will focus ourselves, on the case when $S_\omega$ is replaced by the finite sector:  
\begin{equation*}
S_{\omega,\rho }=\left\{ (x,y)=(r\cos \theta ,r\sin \theta ):0<r<\rho \text{ and } 0<\theta <\omega \right\},
\end{equation*}
and its two lateral edges: 
\begin{equation*}
\left\{ 
\begin{array}{lll}
\Gamma _{0} &=& (0,\rho)  \times \left\{ 0\right\} \\ 
\Gamma _{\omega,\rho } &=& \left\{ (r\cos \omega ,r\sin \omega )~:~0<r<\rho
\right\} ,
\end{array}
\right.
\end{equation*}
where $\rho > 0$ is given small enough.

Then, in polar coordinates, problem \eqref{Pb cone infini} writes 
\begin{equation}\label{pb v}
\left\{\begin{array}{l}
\Lambda_2 v - k \Lambda_1 v = g \quad \text{in } S_{\omega,\rho} \\ \ecart
v(r,0) = v(r,\omega ) = \dfrac{\partial v}{\partial \theta }(r,0) = \dfrac{\partial v}{\partial \theta }(r,\omega) = 0 \\ 
v(\rho,\theta) \text{ given},
\end{array}\right.
\end{equation}
where $\Lambda_1$ and $\Lambda_2$ are given by \eqref{Lambda 1} and \eqref{Lambda 2}.

In this article, we will focus ourselves on the resolution of problem \eqref{pb v} to obtain the behavior of the solution $v$ to problem \eqref{pb v} in $L^p$-weighted spaces, in the neighborhood of the top of the cone. To this $v$ corresponds a solution $u_0$ to problem \eqref{Pb cone infini} by applying the inverse changes of variables and functions. We point out, that there is no reason for this solution $u_0$ to coincide with the variational solution $u_{var}$ belonging to the Hilbert space $H^2(S_{\omega,\rho}) \cap H^1_0(S_{\omega,\rho})$.

To solve problem \eqref{pb v}, we will use results given in  \cite{Cone P2} and to this end, we need to consider 
$$\tau = \min_{j\geqslant 1}\left|\text{Im}(z_j)\right| > 0,$$
where $(z_j)_{j\geqslant 1}$ are the solutions of the following transcendent equation
$$\left(\sinh(z) + z\right) \left(\sinh(z) - z\right) = 0, \quad \text{with} \quad\text{Re}(z)>0.$$
Actually, according to \cite{fadle}, we have
$$\tau \simeq 4.21239.$$
We assume that 
\begin{equation}\label{hyp inv sum 0}
\omega\left(3 - \frac{2}{p}\right) < \tau.
\end{equation}
Recall that $\dis 3-\frac{2}{p}$ is exactly the Sobolev exponent of the space $W^{3,p}$ in two variables. 

Let us remark that we have two cases: 
\begin{enumerate}
\item if $\dis 0 < \omega \leqslant \frac{\tau}{3} \simeq 0.45\pi$, then \eqref{hyp inv sum 0} is satisfied for all $p \in (1,+\infty)$.

\item if $\dis \frac{\tau}{3} < \omega < \tau \simeq 1.34\pi$, then \eqref{hyp inv sum 0} is satisfied for $\dis 1 < p < \frac{2\omega}{3\omega - \tau}$.
\end{enumerate}
Our main result is the following:
\begin{Th}\label{Th principal}
There exists $\rho_0 > 0$ such that for all $\rho \in (0, \rho_0]$, $g \in L^p_{\frac{1}{p}}(S_{\omega,\rho})$ with $p$ satisfying assumption \eqref{hyp inv sum 0}, problem \eqref{pb v} has a unique strong solution 
$$v \in W^{4,p}_{3 - \frac{1}{p}}(S_{\omega,\rho}).$$
More precisely, for $i=0,1,2$ and $j = 0,1,2,3,4$ such that $0\leqslant i+j\leqslant 4$, we have
$$ \frac{\partial^{i+j} v}{\partial r^i \partial \theta^j} \in L^p(S_{\omega,\rho}),$$ 
and
$$\frac{\partial^3 v}{\partial r^3}, \frac{\partial^4 v}{\partial r^3 \partial \theta} \in L^p_{2-\frac{1}{p}}(S_{\omega,\rho}) \quad \text{and} \quad \frac{\partial^4 v}{\partial r^4} \in L^p_{3-\frac{1}{p}}(S_{\omega,\rho}).$$
\end{Th}

\section{Reformulations of problem \eqref{pb v}} \label{Sect Reformulation}

\subsection{Some preliminary calculus}

Let us introduce the auxiliary function $\displaystyle \frac{v}{r}$. Then 
\begin{eqnarray*}
\left( r\dfrac{\partial }{\partial r}\right) ^{2}\left( \frac{v}{r}\right) &=&\left( r\dfrac{\partial }{\partial r}\right) \left( r\dfrac{\partial }{%
\partial r}\right) \left(\frac{v}{r}\right) \\ 
&=&\frac{v}{r}-\dfrac{\partial v}{\partial r}+r\dfrac{\partial ^{2}v}{\partial r^{2}}
\\
&=&\frac{v}{r}-2\dfrac{\partial }{\partial r}\left( r.\frac{v}{r}\right) +r\left( \frac{1}{r}%
\dfrac{\partial v}{\partial r}+\dfrac{\partial ^{2}v}{\partial r^{2}}\right) 
\\
&=& - \frac{v}{r} - 2r\dfrac{\partial }{\partial r}\left( \frac{v}{r}\right) +r\left( \frac{1}{r}%
\dfrac{\partial v}{\partial r}+\dfrac{\partial ^{2}v}{\partial r^{2}}\right),
\end{eqnarray*}
so
\begin{equation*}
\Delta u =\left( \frac{1}{r}\dfrac{\partial v}{\partial r}+\dfrac{\partial^{2} v}{\partial r^{2}}\right) +\frac{1}{r^{2}}\dfrac{\partial ^{2}v}{\partial \theta ^{2}} = \frac{1}{r}\left[ \left( r\dfrac{\partial }{\partial r}\right) ^{2}\left(\frac{v}{r}\right) + \frac{v}{r} + 2r\dfrac{\partial }{\partial r}\left( \frac{v}{r}\right) \right] + \frac{1}{r}\dfrac{\partial ^{2}}{\partial \theta ^{2}}\left( \frac{v}{r}\right).
\end{equation*}
Moreover, we have
\begin{eqnarray*}
\Pi_1 & := &\dis \left( r\dfrac{\partial }{\partial r}\right) ^{2}\left[ \left( r\dfrac{\partial }{\partial r}\right) ^{2}\left[ \frac{v}{r}\right] \right] \\
 & = & \left( r\dfrac{\partial }{\partial r}\right) ^{2}\left[\left( r\dfrac{\partial }{\partial r}\right) 
\left( -\frac{v}{r}+\dfrac{\partial v}{\partial r}\right) \right] \\
 &=& \left( r\dfrac{\partial }{\partial r}\right) ^{2}\left[ \frac{v}{r}-%
\dfrac{\partial v}{\partial r}+r\dfrac{\partial ^{2}v}{\partial r^{2}}\right]
\\
&=&\left( r\dfrac{\partial }{\partial r}\right) \left[ -\frac{v}{r}+\dfrac{%
\partial v}{\partial r}+r^{2}\dfrac{\partial ^{3}v}{\partial r^{3}}\right] 
\\
&=&\frac{v}{r}-\dfrac{\partial v}{\partial r}+r\dfrac{\partial ^{2}v}{%
\partial r^{2}}+2r^{2}\dfrac{\partial ^{3}v}{\partial r^{3}}+r^{3}\dfrac{%
\partial ^{4}v}{\partial r^{4}},
\end{eqnarray*}
and
\begin{eqnarray*}
\Pi_2 & := & 2\left( \dfrac{\partial ^{2}}{\partial \theta ^{2}}-1\right) \left( r \dfrac{\partial }{\partial r}\right) ^{2}\left( \frac{v}{r}\right) +\left( \dfrac{\partial ^{2}}{\partial \theta ^{2}}+1\right) ^{2}\left( \frac{v}{r}\right)  \\ 
&=& 2 \left( \dfrac{\partial ^{2}}{\partial \theta ^{2}}-1\right) \left( \frac{v}{r} -
\dfrac{\partial v}{\partial r}+r\dfrac{\partial ^{2}v}{\partial r^{2}}%
\right) +\frac{1}{r}\left( \dfrac{\partial ^{2}}{\partial \theta ^{2}}%
+1\right) ^{2}v \\
&=&2\dfrac{\partial v}{\partial r}-2r\dfrac{\partial ^{2}v}{\partial r^{2}}+%
\frac{4}{r}\dfrac{\partial ^{2}v}{\partial \theta ^{2}}-2\dfrac{\partial
^{3}v}{\partial \theta ^{2}\partial r}+2r\dfrac{\partial ^{4}v}{\partial
\theta ^{2}\partial r^{2}}+\frac{1}{r}\dfrac{\partial ^{4}v}{\partial \theta
^{4}}-\frac{v}{r}.
\end{eqnarray*}
Therefore, we obtain that
$$r^3 \Delta^2 u = \Pi_1 + \Pi_2.$$
We set 
\begin{equation*}
w=\frac{v}{r}.
\end{equation*}
Then, in $S_{\omega,\rho}$, problem \eqref{pb v} writes 
\begin{equation} \label{EquationEn RetTheta}
\left\{ 
\begin{array}{l}
\dfrac{1}{r^{3}}\left[ \left( r\dfrac{\partial }{\partial r}\right) ^{2}%
\left[ \left( r\dfrac{\partial }{\partial r}\right) ^{2}w\right] +2\left( 
\dfrac{\partial ^{2}}{\partial \theta ^{2}}-1\right) \left( r\dfrac{\partial 
}{\partial r}\right) ^{2}w+\left( \dfrac{\partial ^{2}}{\partial \theta ^{2}}%
+1\right) ^{2}w\right] \\ \ecart
-\dfrac{k}{r}\left[ \left( \left( r\dfrac{\partial }{\partial r}\right)
^{2}w+2\left( r\dfrac{\partial }{\partial r}\right) w+w\right) +\dfrac{%
\partial ^{2}w}{\partial \theta ^{2}}\right] =g \\ \\
w(r,0)=w(r,\omega )=\dfrac{\partial w}{\partial \theta }(r,0)=\dfrac{%
\partial w}{\partial \theta }(r,\omega )=0 \\
w(\rho, \theta) \text{ given}. 
\end{array}\right. 
\end{equation}

\subsection{Vector formulation of problem \eqref{EquationEn RetTheta}}

Now, let us consider the vector variable $\Psi (r,\theta ):$ 
\begin{equation*}
\Psi =\left( 
\begin{array}{c}
w \\ \ecart
\displaystyle \left( r\dfrac{\partial }{\partial r}\right) ^{2}w
\end{array}%
\right) ,
\end{equation*}%
and the following matrix  
\begin{equation*}
\mathcal{A}=\left( 
\begin{array}{cc}
0 & 1 \\ \ecart
-\left( \dfrac{\partial ^{2}}{\partial \theta ^{2}}+1\right) ^{2} & -2\left( 
\dfrac{\partial ^{2}}{\partial \theta ^{2}}-1\right) 
\end{array}%
\right).
\end{equation*}
We have
\begin{eqnarray*}
&&\hspace*{-1cm}\left( r\dfrac{\partial }{\partial r}\right) ^{2}\Psi -\mathcal{A}\Psi  \\
&=&\left( 
\begin{array}{c}
\displaystyle \left( r\dfrac{\partial }{\partial r}\right) ^{2}w  \\ \ecart 
\displaystyle \left( r\dfrac{\partial }{\partial r}\right) ^{2}\left( r\dfrac{\partial }{
\partial r}\right) ^{2}w 
\end{array}
\right) +\left( 
\begin{array}{c}
\displaystyle -\left( r\dfrac{\partial }{\partial r}\right) ^{2} w  \\ \ecart
\displaystyle \left( \dfrac{\partial ^{2}}{\partial \theta ^{2}}+1\right) ^{2} w + 2 \left( \dfrac{\partial ^{2}}{\partial \theta ^{2}}-1\right)
\left( r\dfrac{\partial }{\partial r}\right) ^{2}w 
\end{array}%
\right)  \\ \ecart
&=&\left( 
\begin{array}{c}
0 \\ \ecart
\displaystyle \left( r\dfrac{\partial }{\partial r}\right) ^{2}\left( r\dfrac{\partial }{\partial r}\right) ^{2}w +2\left( \dfrac{\partial ^{2}}{\partial \theta ^{2}}-1\right) \left( r\dfrac{\partial }{\partial r}\right)^{2}w +\left( \dfrac{\partial ^{2}}{\partial \theta ^{2}}+1\right) ^{2}w 
\end{array}
\right) \\ \ecart
& = & \binom{0}{\Pi_1 +\Pi_2 } = \binom{0}{r^3 \Delta^2 u}.
\end{eqnarray*}

We set
$$\A_0 = \left( 
\begin{array}{cc}
0 & 0 \\ 
\dfrac{\partial ^{2}}{\partial \theta ^{2}}+1 & 1%
\end{array}%
\right) \quad \text{and}\quad \B_0 = \left( 
\begin{array}{cc}
0 & 0 \\ 
2\left( r\dfrac{\partial }{\partial r}\right)  & 0%
\end{array}\right).$$
We will precise the domain of all these operators in section \ref{sect espace de travail et domaine}. It is clear that the action of these operators are independent. Then
$$
\left(\mathcal{A}_0 + \mathcal{B}_0\right)\Psi = \left( 
\begin{array}{cc}
0 & 0 \\ \ecart
\displaystyle\dfrac{\partial ^{2}}{\partial \theta ^{2}}+1 & 1%
\end{array}%
\right) \left( 
\begin{array}{c}
w \\ \ecart
\displaystyle\left( r\dfrac{\partial }{\partial r}\right) ^{2} w
\end{array}%
\right) + \left( 
\begin{array}{cc}
0 & 0 \\ \ecart
2\left( r\dfrac{\partial }{\partial r}\right)  & 0%
\end{array}%
\right)  \left( 
\begin{array}{c}
w \\ \ecart
\displaystyle\left( r\dfrac{\partial }{\partial r}\right) ^{2}w 
\end{array}%
\right)$$ 
hence
\begin{eqnarray*}
\left(\mathcal{A}_0 + \mathcal{B}_0\right) \Psi &=& \left( 
\begin{array}{c}
0 \\ \ecart
\displaystyle\dfrac{\partial ^{2}}{\partial \theta ^{2}} w + w +\left( r \dfrac{\partial }{\partial r}\right) ^{2}w 
\end{array}
\right) + \left( 
\begin{array}{c}
0 \\ \ecart
\displaystyle2\left( r\dfrac{\partial }{\partial r}\right) w 
\end{array}
\right) \\ \ecart
&=&\left( 
\begin{array}{c}
0 \\ \ecart
\displaystyle 2r\dfrac{\partial }{\partial r} w + \dfrac{\partial ^{2}}{%
\partial \theta ^{2}}w + w + \left( r\dfrac{\partial }{\partial r}\right) ^{2} w 
\end{array}
\right)  \\ \ecart
&=&\left( 
\begin{array}{c}
0 \\ 
r\Delta u
\end{array}
\right) = r \left( 
\begin{array}{c}
0 \\ 
\Delta u
\end{array}
\right) .
\end{eqnarray*}
The generalized diffusion equation writes as  
\begin{equation*}
\dfrac{1}{r^{3}}\left[ \left( r\dfrac{\partial }{\partial r}\right) ^{2}\Psi
-\mathcal{A}\Psi \right] -\dfrac{k}{r}\left(\mathcal{A}_0 + \mathcal{B}_0\right)\Psi =\left( 
\begin{array}{c}
0 \\ 
g
\end{array}
\right) .
\end{equation*}
Finally, we obtain the following complete equation 
\begin{equation}\label{eq complete}
\left[ \left( r\dfrac{\partial }{\partial r}\right) ^{2}\Psi -\mathcal{A}%
\Psi \right] -kr^{2}\mathcal{A}_{0}\Psi -kr^{2}\mathcal{B}_{0}\Psi =\left( 
\begin{array}{c}
0 \\ 
r^{3}g%
\end{array}%
\right) .
\end{equation}
Note that linear operators $\mathcal{A}$ and $\mathcal{A}_{0}$ act with respect to variable $\theta $ whereas operator $\mathcal{B}_{0}$ acts with respect to variable $r\dfrac{\partial }{\partial r}$.

\subsection{New formulation in a finite cone}

We apply the following variables and functions change  
\begin{equation*}
r = \rho e^{-t}, \quad \phi (t,\theta )= w(\rho
e^{-t},\theta ) \quad \text{and}\quad g(\rho e^{-t},\theta
)=G(t,\theta )
\end{equation*}
then, it is easy to verify that 
\begin{equation*}
\left( r\dfrac{\partial }{\partial r}\right) w = - \frac{\partial \phi}{\partial t}  \quad \text{and} \quad \left( r\dfrac{\partial }{\partial r}\right) ^{2} w =\left[r\dfrac{\partial}{\partial r} + r^{2}\dfrac{\partial ^{2}}{\partial r^{2}}\right] w = \frac{\partial^2 \phi}{\partial t^2}.
\end{equation*}
We set $\Phi(t,\theta) = \Psi(\rho e^{-t}, \theta)$; then
\begin{equation*}
\Phi = \left( 
\begin{array}{c}
\phi \\ \ecart
\dis\frac{\partial^2 \phi}{\partial t^2}
\end{array}
\right).
\end{equation*}
Therefore, equation \eqref{eq complete} is now set on the strip $\Sigma = (0,+\infty) \times (0,\omega)$ and writes
\begin{eqnarray*}
&&\hspace*{-2.5cm}\left[ \left( r\dfrac{\partial }{\partial r}\right) ^{2}\Psi -\mathcal{A}
\Psi \right] -kr^{2}\mathcal{A}_{0}\Psi -kr^{2}\mathcal{B}_{0}\left[ \Psi %
\right] \\
&=&\left[ \frac{\partial^2 \Phi}{\partial t^2} -\mathcal{A}\Phi \right] -k\rho ^{2}e^{-2t}\mathcal{A%
}_{0}\Phi +k\rho ^{2}e^{-2t}\mathcal{B}_{1}\left[ \Phi \right] =\left( 
\begin{array}{c}
0 \\ 
\rho ^{3}e^{-3t}G
\end{array}
\right),
\end{eqnarray*}
where
\begin{equation*}
\mathcal{B}_{1}=\left( 
\begin{array}{cc}
0 & 0 \\ 
\displaystyle -2\frac{\partial}{\partial t} & 0%
\end{array}%
\right) .
\end{equation*}
The boundary conditions, in problem \eqref{EquationEn RetTheta}, on $w$ become  
\begin{equation*}
\phi (r,0)=\phi (r,\omega )=\dfrac{\partial \phi }{\partial \theta }(r,0)=%
\dfrac{\partial \phi }{\partial \theta }(r,\omega )=0.
\end{equation*}
As usual, we will use the vector valued notation:
\begin{equation*}
\Phi (t)(\theta) := \Phi (t,\theta )=\left( 
\begin{array}{c}
\phi (t,\theta ) \\ \ecart 
\displaystyle \frac{\partial^2 \phi}{\partial t^2} (t,\theta )%
\end{array}%
\right) =\left( 
\begin{array}{c}
\phi (t,.) \\ \ecart
\displaystyle\frac{\partial^2 \phi}{\partial t^2} (t,.)%
\end{array}%
\right) (\theta ) :=\left( 
\begin{array}{c}
\phi (t) \\ \ecart
\displaystyle\frac{\partial^2 \phi}{\partial t^2} (t)%
\end{array}%
\right) (\theta ).
\end{equation*}
Hence, we deduce the following abstract vector valued equation 
\begin{equation*}
\left[ \Phi'' (t)-\mathcal{A}\Phi (t)\right] -k\rho
^{2}e^{-2t}\mathcal{A}_{0}\Phi (t)+k\rho ^{2}e^{-2t}\left[ \mathcal{B}%
_{1}\Phi \right] (t)=\left( 
\begin{array}{c}
0 \\ 
\rho ^{3}e^{-3t}G(t)%
\end{array}%
\right) ,
\end{equation*}
where
\begin{equation*}
\left[ \mathcal{B}_{1}\Phi \right] (t) = \left( 
\begin{array}{cc}
0 & 0 \\ \ecart
\displaystyle -2 \frac{\partial}{\partial t} & 0%
\end{array}%
\right) \left( 
\begin{array}{c}
\phi \\ \ecart
\displaystyle\frac{\partial^2 \phi}{\partial t^2}
\end{array}
\right) (t) = \left( 
\begin{array}{c}
0 \\ \ecart
\displaystyle -2 \frac{\partial \phi}{\partial t}(t)
\end{array}
\right).
\end{equation*}
Note that $\mathcal{A}$ and $\mathcal{A}_{0}$ act on $\Phi (t)$ with respect to $\theta$, while $\mathcal{B}_{1}$ acts on $\Phi $ with respect to $t$.

To determine completely $\Psi $ in a suitable space and also $\Phi $, it is necessary to give us a boundary condition on $\Psi$ in $\rho$. It means that $\Phi$ is given at $t=0$ : 
\begin{equation*}
\Phi (0)=\Phi _{0}.
\end{equation*}
Now, we have to solve on $(0,+\infty)$ the following problem 
\begin{equation}
\left\{ 
\begin{array}{l}
\Phi ^{\prime \prime }(t)-\mathcal{A}\Phi (t)-k\rho ^{2}e^{-2t}\mathcal{A}%
_{0}\Phi (t)+ k\rho ^{2}e^{-2t} \left[ \mathcal{B}_{1}\Phi \right] (t)=\left( 
\begin{array}{c}
0 \\ 
\rho ^{3}e^{-3t}G(t)%
\end{array}%
\right) \\ 
\Phi (0)=\Phi _{0}.%
\end{array}%
\right.  \label{ProbConeFini}
\end{equation}
\begin{Rem}\label{Rem Phi(+infini)=0}
Note that the boundary condition at $t=+\infty$, will be included in the vector valued space containing the solution $\Phi$.
\end{Rem}

\subsection{Sums of linear operators} \label{sect espace de travail et domaine}

In this section, we are going to write problem \eqref{ProbConeFini} as a sum of linear operators, firstly in the following Banach space 
\begin{equation*}
X=W_{0}^{2,p}(0,\omega) \times L^{p}(0,\omega),
\end{equation*}
see \eqref{ProbAbstrait} below and secondly, in $L^p(0,+\infty;X)$, see \eqref{Pb L1+L2} below.

Here $X$ is equipped, for instance, with the following norm 
\begin{equation*}
\left\Vert \left( 
\begin{array}{c}
\psi _{1} \\ 
\psi _{2}%
\end{array}%
\right) \right\Vert _{X}=\left\Vert \psi _{1}\right\Vert
_{W_{0}^{2,p}(0,\omega)}+\left\Vert \psi _{2}\right\Vert
_{L^{p}(0,\omega)},
\end{equation*}
where
\begin{equation*}
W_{0}^{2,p}(0,\omega) = \left\{ \varphi \in W^{2,p}(0,\omega) : \varphi (0) = \varphi (\omega ) = \varphi'(0) = \varphi'(\omega) = 0\right\} .
\end{equation*}
Then, we define the linear operator $\mathcal{A}$ by
\begin{equation*}
\left\{ 
\begin{array}{lll}
D(\mathcal{A}) & = & \dis \left[ W^{4,p}(0,\omega)\cap W_{0}^{2,p}(0,\omega)\right] \times W_{0}^{2,p}(0,\omega)\subset X \\ \ecart
\mathcal{A}\left( 
\begin{array}{c}
\psi _{1} \\ 
\psi _{2}
\end{array}
\right) &=& \left( 
\begin{array}{c}
\psi _{2} \\ 
-\left( \dfrac{\partial ^{2}}{\partial \theta ^{2}}+1\right) ^{2}\psi
_{1}-2\left( \dfrac{\partial ^{2}}{\partial \theta ^{2}}-1\right) \psi _{2}%
\end{array}%
\right), \quad \left( 
\begin{array}{c}
\psi _{1} \\ 
\psi _{2}%
\end{array}%
\right) \in D(\A).
\end{array}
\right.
\end{equation*}
In the same way, we define operator $\mathcal{A}_{0}$ by 
\begin{equation*}
\left\{ 
\begin{array}{lll}
D(\mathcal{A}_{0}) & = & W_{0}^{2,p}(0,\omega) \times L^{p}(0,\omega) = X \\ \ecart
\mathcal{A}_{0}\left( 
\begin{array}{c}
\psi _{1} \\ 
\psi _{2}
\end{array}
\right) & = & \dis \left( 
\begin{array}{c}
0 \\ 
\left( \dfrac{\partial ^{2}}{\partial \theta ^{2}}+1\right) \psi _{1}+\psi
_{2}%
\end{array}%
\right), \quad \left( 
\begin{array}{c}
\psi _{1} \\ 
\psi _{2}%
\end{array}%
\right) \in D(\A_0).%
\end{array}%
\right.
\end{equation*}
It is clear that $D(\mathcal{A})\subset D(\mathcal{A}_{0})$. Note that operator $\mathcal{A}_{0}$ is continuous from $X$ into itself since
\begin{eqnarray*}
\left\Vert \mathcal{A}_{0}\left( 
\begin{array}{c}
\psi _{1} \\ 
\psi _{2}%
\end{array}%
\right) \right\Vert _{X} &=&\left\Vert \left( 
\begin{array}{c}
0 \\ 
\left( \dfrac{\partial ^{2}}{\partial \theta ^{2}}+1\right) \psi _{1}+\psi
_{2}
\end{array}
\right) \right\Vert _{X}=\left\Vert \left( \dfrac{\partial ^{2}}{\partial
\theta ^{2}}+1\right) \psi _{1}+\psi _{2}\right\Vert _{L^{p}(0,\omega)} \\
&\leqslant &\left\Vert \psi _{1}\right\Vert _{W_{0}^{2,p}(0,\omega)} + \left\Vert \psi _{2}\right\Vert _{L^{p}(0,\omega)}=\left\Vert
\left( 
\begin{array}{c}
\psi _{1} \\ 
\psi _{2}%
\end{array}%
\right) \right\Vert _{X}.
\end{eqnarray*}
Equation \eqref{ProbConeFini} is set in the Banach space $X$.

Recall that the second member in problem \eqref{Pb cone infini} satisfies
$$f \in L^p\left(S_{\omega,\rho}\right), \quad \text{for } p \in (1,+\infty).$$
Set
\begin{equation*}
t\mapsto e^{-3t}G(t)(.)=e^{-3t}G(t,.)=H(t,.)=H(t)(.).
\end{equation*}
Therefore, we have 
\begin{eqnarray*}
\int_{S_{\omega,\rho }}\left\vert f(x,y)\right\vert ^{p}dxdy &=&\int_{S_{\omega,\rho
}}\left\vert g(r,\theta )\right\vert ^{p}rdrd\theta =\rho ^{2}\int_{%
\Sigma }\left\vert G(t,\theta )\right\vert ^{p}e^{-2t}dtd\theta \\
&=&\rho ^{2}\int_{\Sigma }\left\vert e^{\left( 3-\frac{2}{p}\right) t}H(t,\theta )\right\vert ^{p}dtd\theta \\
&=&\rho ^{2}\int_{0}^{+\infty }\left\vert e^{\left( 3-\frac{2}{p}\right)
t}\right\vert ^{p}\left[ \left( \int_{0}^{\omega }\left\vert H(t)(\theta
)\right\vert ^{p}d\theta \right) ^{1/p}\right] ^{p}dt \\
&=&\rho ^{2}\int_{0}^{+\infty }\left[ e^{\left( 3-\frac{2}{p}\right)
t}\left\Vert H(t)\right\Vert _{L^{p}(0,\omega)}\right] ^{p}dt.
\end{eqnarray*}
It follows that the second member $H$ is in the weighted Sobolev space
\begin{equation*}
\left\{H : t\mapsto e^{\left( 3-\frac{2}{p}\right) t}H\in L^{p}(\Sigma)\right\} = L_{\nu}^{p}(0,+\infty ;L^{p}(0,\omega)).
\end{equation*}
where 
\begin{equation*}
\nu = 3-\frac{2}{p} > 1,
\end{equation*}
is exactly the Sobolev exponent of the space $W^{3,p}(\Sigma)$ in dimension 2.

Then, since it would not be easy to work in weighted Sobolev spaces, we will use the following new vector valued function :
\begin{equation}\label{V = Phi}
V(t)=e^{\nu t}\Phi (t) = \left( 
\begin{array}{c}
e^{\nu t}\phi (t) \\ 
e^{\nu t} \phi''(t)
\end{array}%
\right) =\left( 
\begin{array}{c}
V_{1}(t) \\ 
V_{2}(t)
\end{array}
\right) .
\end{equation}
Since we have
\begin{equation*}
\Phi (t)=\left( 
\begin{array}{c}
e^{-\nu t}V_{1}(t) \\ 
e^{-\nu t}V_{2}(t)
\end{array}
\right),
\end{equation*}
we deduce that
\begin{equation*}
\Phi'(t) = -\nu e^{-\nu t}V(t)+e^{-\nu t}V'(t) = e^{-\nu t}\left(
\partial_{t}-\nu I\right) V(t) ,
\end{equation*}
and
\begin{equation*}
\Phi''(t) = \nu ^{2}e^{-\nu t}V(t) - 2\nu e^{-\nu t} V'(t) + e^{-\nu t} V''(t) = e^{-\nu t}\left( \partial_{t}-\nu I\right)^{2}V(t).
\end{equation*}
Moreover, we obtain
\begin{eqnarray*}
\left[ \mathcal{B}_{1}\Phi \right] (t) &=&\left( 
\begin{array}{cc}
0 & 0 \\ 
-2\partial_{t} & 0
\end{array}\right) \left( 
\begin{array}{c}
e^{-\nu t}V_{1}(t) \\ 
e^{-\nu t}V_{2}(t)%
\end{array}%
\right) \\
&=&e^{-\nu t}\left( 
\begin{array}{cc}
0 & 0 \\ 
-2(\partial_{t}-\nu I) & 0%
\end{array}%
\right) \left( 
\begin{array}{c}
V_{1}(t) \\ 
V_{2}(t)%
\end{array}%
\right) \\
&=&e^{-\nu t}\left[ \mathcal{B}_{2}V\right] (t)
\end{eqnarray*}%
where
\begin{equation*}
\mathcal{B}_{2}=\left( 
\begin{array}{cc}
0 & 0 \\ 
-2(\partial_{t}-\nu I) & 0%
\end{array}%
\right) .
\end{equation*}
Hence, problem \eqref{ProbConeFini} becomes
\begin{equation*}
\left\{ 
\begin{array}{l}
e^{-\nu t}\left( \partial_{t}-\nu I\right)^{2}V(t)-e^{-\nu t}\mathcal{A} V(t)-k\rho ^{2}e^{-\nu t}e^{-2t}\mathcal{A}_{0}V(t) \\ 
-k\rho ^{2}e^{-\nu t}e^{-2t}\left[ (\mathcal{B}_{2}V)\right] (t)=\left( 
\begin{array}{c}
0 \\ 
\rho ^{3}H(t)
\end{array}
\right) \\ 
V(0)=\Phi _{0},
\end{array}
\right.
\end{equation*}
then
\begin{equation*}
\left\{ 
\begin{array}{l}
\left( \partial_{t}-\nu I\right) ^{2}V(t)-\mathcal{A}V(t)-k\rho ^{2}e^{-2t}
\mathcal{A}_{0}V(t) -k\rho ^{2}e^{-2t}\left[ (\mathcal{B}_{2}V)\right] (t)=\left( 
\begin{array}{c}
0 \\ 
\rho ^{3}e^{\nu t}H(t)
\end{array}\right) \\ 
V(0)=\Phi _{0}.
\end{array}
\right.  
\end{equation*}
Note that \refR{Rem Phi(+infini)=0} holds true for $V$.
 
We will be interested, in particular, in the following homogeneous problem
\begin{equation}\label{ProbAbstrait}
\left\{ \hspace*{-0.1cm}
\begin{array}{l}
\left( \partial_{t}-\nu I\right) ^{2}V(t)-\mathcal{A}V(t)-k\rho ^{2}e^{-2t}
\mathcal{A}_{0}V(t) -k\rho ^{2}e^{-2t}\left[ (\mathcal{B}_{2}V)\right] (t)=\left( 
\begin{array}{c}
0 \\ 
\rho ^{3}e^{\nu t}H(t)
\end{array}\right) \\ 
V(0)=0.
\end{array}
\right.  
\end{equation}
We have 
\begin{equation*}
t\mapsto \rho ^{3}e^{\nu t}H(t)\in L^{p}(\Sigma) = L^{p}(0,+\infty ;L^{p}(0,\omega)),
\end{equation*}
and
\begin{equation*}
t\mapsto \left( 
\begin{array}{c}
0 \\ 
\rho ^{3}e^{\nu t}H(t)
\end{array}
\right) \in L^{p}\left( 0,+\infty ;W_{0}^{2,p}(0,\omega)\times L^{p}(0,\omega)\right) = L^{p}(0,+\infty ;X).
\end{equation*}
Finally, let us introduce the following abstract linear operators: 
\begin{equation*}
\left\{ \begin{array}{cll}
D(\mathcal{L}_{1}) & = & \dis \left\{ V\in W^{2,p}(0,+\infty ;X): V(0) = V(+\infty) = 0\right\} \\ \ecart
\left[ \mathcal{L}_{1}(V)\right] (t) & = & \dis \left( \partial_{t}-\nu I\right)^{2}V(t)=V''(t)-2\nu V'(t)+\nu^{2}V(t),
\end{array}\right.
\end{equation*}
with $\nu = 3 - \frac{2}{p} \in (1,3)$, 
\begin{equation*}
\left\{ 
\begin{array}{lll}
D(\mathcal{L}_{2}) & = &\dis \left\{ V\in L^{p}(0,+\infty ;X): \text{for }a.e.~t\in
(0,+\infty ),~ V(t)\in D(\mathcal{A})\right\} \\ \ecart
\left[ \mathcal{L}_{2}(V)\right] (t) & = & -\mathcal{A}V(t),
\end{array}
\right.
\end{equation*}
\begin{equation*}
\left\{ 
\begin{array}{lll}
D(\mathcal{P}_{1}) & = &\dis \left\{ V\in L^{p}(0,+\infty ;X): \text{for }a.e.~t\in
(0,+\infty ),~ V(t)\in D(\mathcal{A}_0)\right\} \\ \ecart
\left[ \mathcal{P}_{1}(V)\right] (t) & = & -e^{-2t} \mathcal{A}_{0}V(t),
\end{array}
\right.
\end{equation*}
and
\begin{equation*}
\left\{ 
\begin{array}{lll}
D(\mathcal{P}_{2}) & = &\dis W^{1,p}(0,+\infty ;X) \\ \ecart
\left[ \mathcal{P}_{2}(V)\right] (t) & = & -e^{-2t}\left( \mathcal{B}_{2}V\right) (t).
\end{array}\right.
\end{equation*}
Then, problem \eqref{ProbAbstrait} can be written as the following abstract equation
\begin{equation}\label{Pb L1+L2}
\left(\L_1 + \L_2\right) V + k \rho^2 \left(\P_1 + \P_2 \right)V = \F,
\end{equation}
set in $L^p(0,+\infty;X)$, with $p \in (1,+\infty)$, where, for almost every $t \in (0,+\infty)$
$$\F(t) =  \left( 
\begin{array}{c}
0 \\ 
\rho ^{3}e^{\sigma t}H(t)
\end{array}\right).$$

\section{Proof of \refT{Th principal}} \label{Sect proof of main Th}

\subsection{Resolution of equation \eqref{Pb L1+L2}}

Equation \eqref{Pb L1+L2} will be completely studied in the second part of this work by using the sum theory of linear operators, where the main result described by Theorem 1.1 in \cite{Cone P2} states that there exists $\rho_0 >0$ such that for all $\rho \in(0,\rho_0]$, there exists a unique solution $V \in D(\L_1+\L_2)$ to equation \eqref{Pb L1+L2} that is
$$V \in W^{2,p}(0,+\infty;X) \cap L^p(0,+\infty;D(\A)).$$
Thus, we know that there exists a continuous extension from $W^{2,p}(0,+\infty;X)$ into $W^{2,p}(\RR;X)$ and also from $L^p(0,+\infty;D(\A))$ into $L^p(\RR;D(\A))$; it suffices, for instance to use the well-known Babich techniques.

Set $\widetilde{V}$, this extension of $V$ which writes
$$\widetilde{V} = \left(\begin{array}{c}
\widetilde{V_1} \\
\widetilde{V_2}
\end{array}\right).$$
Then,  
$$\widetilde{V_1} \in W^{2,p}\left(\RR; W^{2,p}_0(0,\omega)\right) \cap L^p\left(\RR; W^{4,p}(0,\omega)\cap  W^{2,p}_0(0,\omega)\right),$$
and
$$\widetilde{V_2} \in W^{2,p}\left(\RR; L^p(0,\omega)\right) \cap L^p\left(\RR; W^{2,p}_0(0,\omega)\right).$$
We deduce, for instance, that 
$$\widetilde{V_1}, \widetilde{V_2} \in W^{2,p}\left(\RR \times (0,\omega)\right) ,$$
by using the Mihlin's theorem (see \cite{mihlin}) and 
$$\widetilde{V_1} \in W^{2,p}\left(\RR; W^{2,p}_0(0,\omega)\right) \cap L^p\left(\RR; W^{4,p}(0,\omega)\right).$$
Therefore, we deduce that $V_1$ and $V_2$ have the same regularities on $(0,\rho) \times (0,\omega)$, with $\rho \in (0,\rho_0]$.
 
\subsection{Regularity of $V(t,\theta)=V(t)(\theta)$}

Recall that, from \eqref{V = Phi}, we have
$$V(t)=e^{\nu t}\Phi (t) \quad \text{and} \quad \Phi (t) = e^{-\nu t}V(t).$$
Moreover, since $r = \rho e^{-t}$, we have
$$V_1 (t,\theta) = \frac{v(\rho e^{-t},\theta)}{\left(\rho e^{-t}\right)^{\nu+1}},$$
where $\dis\nu = 3 - \frac{2}{p} \in (1,3)$ and
$$V_2 (t,\theta) = \left(r \frac{\partial}{\partial r}\right)^2 \left(\frac{v}{r}\right)(\rho e^{-t},\theta).$$

\subsubsection{Regularity of $V_1$}\label{SubSect First situation}

Here, we explicit the fact that 
\begin{equation}\label{V1 in W2p}
V_1 \in W^{2,p}\left((0,+\infty)\times(0,\omega)\right).
\end{equation}
We have
$$\int_0^{+\infty} \int_0^\omega \left|V_1(t,\theta)\right|^p d\theta~dt = \int_0^{+\infty} \int_0^\omega \left|\frac{v(\rho e^{-t},\theta)}{\rho^{\nu+1} e^{-t(\nu+1)}}\right|^p d\theta~dt < + \infty.$$
Setting $r = \rho e^{-t}$, we obtain
$$\begin{array}{lll}
\dis\int_0^{+\infty} \int_0^\omega \left|\frac{v(\rho e^{-t},\theta)}{\left(\rho e^{-t}\right)^{\nu+1}}\right|^p d\theta~dt & = & \dis\int_0^\rho \int_0^\omega \left|\frac{v\left(r,\theta\right)}{r^{\nu + 1 + \frac{1}{p}}}\right|^p d\theta~dr \\ \ecart

& = & \dis \int_0^\rho \int_0^\omega r^{-4p + 1}\left|v\left(r,\theta\right)\right|^p d\theta~dr.
\end{array}$$
Then, we have
\begin{equation}\label{v in Lp gamma1}
v \in L^p_{\gamma_0}(S_{\omega,\rho}),
\end{equation}
where $\gamma_0 = -4 + \frac{1}{p}$. Moreover
$$\begin{array}{lll}
\dis\frac{\partial V_1}{\partial t}(t,\theta) &=& \dis \frac{1}{\rho^{\nu + 1}}\frac{\partial}{\partial t}\left(v(\rho e^{-t},\theta) e^{(\nu+1)t}\right) \\ \ecart
& = & \dis\frac{1}{\rho^{\nu + 1}}\left((\nu+1)e^{(\nu+1)t} v(\rho e^{-t},\theta) - \rho e^{\nu t} \frac{\partial v}{\partial r}(\rho e^{-t},\theta)  \right),
\end{array}$$
hence
$$ \frac{1}{\left(\rho e^{- t}\right)^{\nu}} \frac{\partial v}{\partial r}(\rho e^{-t},\theta) = \frac{(\nu+1)}{\left(\rho e^{-t}\right)^{\nu+1}} v(\rho e^{-t},\theta) - \frac{\partial V_1(t,\theta)}{\partial t}.$$
Thus, in virtue of \eqref{V1 in W2p} and \eqref{v in Lp gamma1}, it follows that
$$(t,\theta) \longmapsto \frac{1}{\left(\rho e^{- t}\right)^{\nu}} \frac{\partial v}{\partial r}(\rho e^{-t},\theta) \in L^p((0,+\infty)\times (0,\omega)),$$
and
$$\begin{array}{lll}
\dis\int_0^{+\infty} \int_0^\omega \left|\frac{1}{\left(\rho e^{- t}\right)^{\nu}} \frac{\partial v}{\partial r}(\rho e^{-t},\theta)\right|^p d\theta~dt & = & \dis \int_0^{\rho} \int_0^\omega \left|\frac{1}{r^{\nu + \frac{1}{p}}} \frac{\partial v}{\partial r}(r,\theta)\right|^p d\theta~dr \\ \ecart

& = & \dis \int_0^{\rho} \int_0^\omega r^{-3p + 1} \left| \frac{\partial v}{\partial r}(r,\theta)\right|^p d\theta~dr.
\end{array}$$
So, we obtain
\begin{equation}\label{v' in Lp gamma2}
\frac{\partial v}{\partial r} \in L^p_{\gamma_1}(S_{\omega,\rho}),
\end{equation}
where $\gamma_1 = -3 + \frac{1}{p}$. Furthermore, we have
$$\begin{array}{lll}
\dis\frac{\partial^2 V_1}{\partial t^2}(t,\theta) &=& \dis \frac{\nu+1}{\rho^{\nu + 1}}\frac{\partial}{\partial t}\left(v(\rho e^{-t},\theta) e^{(\nu+1)t}\right) - \frac{\nu}{\rho^{\nu}} e^{\nu t} \frac{\partial v}{\partial r}(\rho e^{-t},\theta) + \frac{\rho}{\rho^{\nu}} e^{-t}e^{\nu t} \frac{\partial^2 v}{\partial r^2}(\rho e^{-t},\theta)  \\ \\
& = & \dis\frac{\nu+1}{\rho^{\nu + 1}}\left((\nu+1)e^{(\nu+1)t} v(\rho e^{-t},\theta) - \rho e^{\nu t} \frac{\partial v}{\partial r}(\rho e^{-t},\theta)  \right) - \frac{\nu}{\rho^{\nu}} e^{\nu t} \frac{\partial v}{\partial r}(\rho e^{-t},\theta) \\ \ecart
&&\dis + \frac{1}{\rho^{\nu -1}} e^{(\nu - 1)t} \frac{\partial^2 v}{\partial r^2}(\rho e^{-t},\theta) \\ \\

& = & \dis \frac{(\nu+1)^2}{\left(\rho e^{-t}\right)^{\nu + 1}}  v(\rho e^{-t},\theta) - \frac{2\nu + 1}{\left(\rho e^{- t}\right)^{\nu}}  \frac{\partial v}{\partial r}(\rho e^{-t},\theta) + \frac{1}{\left(\rho e^{-t}\right)^{\nu-1}}  \frac{\partial^2 v}{\partial r^2}(\rho e^{-t},\theta),
\end{array}$$
hence
$$\frac{1}{\left(\rho e^{-t}\right)^{\nu-1}}  \frac{\partial^2 v}{\partial r^2}(\rho e^{-t},\theta) = \frac{\partial^2 V_1}{\partial t^2}(t,\theta) + \frac{2\nu + 1}{\left(\rho e^{- t}\right)^{\nu}}  \frac{\partial v}{\partial r}(\rho e^{-t},\theta) - \frac{(\nu+1)^2}{\left(\rho e^{-t}\right)^{\nu + 1}}  v(\rho e^{-t},\theta).$$
Thus, in virtue of \eqref{V1 in W2p}, \eqref{v in Lp gamma1} and \eqref{v' in Lp gamma2}, it follows that
$$(t,\theta) \longmapsto \frac{1}{\left(\rho e^{- t}\right)^{\nu - 1}} \frac{\partial^2 v}{\partial r^2}(\rho e^{-t},\theta) \in L^p((0,+\infty)\times (0,\omega)),$$
and
$$\begin{array}{lll}
\dis\int_0^{+\infty} \int_0^\omega \left|\frac{1}{\left(\rho e^{- t}\right)^{\nu -1}} \frac{\partial^2 v}{\partial r^2}(\rho e^{-t},\theta)\right|^p d\theta~dt & = & \dis \int_0^{\rho} \int_0^\omega \left|\frac{1}{r^{\nu - 1 + \frac{1}{p}}} \frac{\partial^2 v}{\partial r^2}(r,\theta)\right|^p d\theta~dr \\ \ecart

& = & \dis \int_0^{\rho} \int_0^\omega r^{-2p + 1} \left| \frac{\partial^2 v}{\partial r^2}(r,\theta)\right|^p d\theta~dr.
\end{array}$$
So, we obtain
\begin{equation*}
\frac{\partial^2 v}{\partial r^2} \in L^p_{\gamma_2}(S_{\omega,\rho}),
\end{equation*}
where $\gamma_2 = -2 + \frac{1}{p}$.

On the other hand
$$\frac{\partial V_1}{\partial \theta} (\rho e^{-t}, \theta) = \frac{1}{\left(\rho e^{-t} \right)^{\nu + 1}}\frac{\partial v}{\partial \theta} (\rho e^{-t}, \theta).$$
Thus, it follows that
$$\begin{array}{lll}
\dis\int_0^{+\infty} \int_0^\omega \left| \frac{\partial V_1}{\partial \theta}(\rho e^{-t},\theta)\right|^p d\theta~dt & = & \dis \int_0^{+\infty} \int_0^\omega \frac{1}{\left(\rho e^{- t}\right)^{p(\nu + 1)}} \left|\frac{\partial v}{\partial \theta}(\rho e^{-t},\theta)\right|^p d\theta~dt \\ \ecart

& = & \dis \int_0^{\rho} \int_0^\omega r^{-4p+1} \left|\frac{\partial v}{\partial \theta}(r,\theta)\right|^p d\theta~dr.
\end{array}$$
So, we obtain
\begin{equation*}
\frac{\partial v}{\partial \theta} \in L^p_{\gamma_0}(S_{\omega,\rho}).
\end{equation*}
In the same way, we deduce that
$$\frac{\partial^2 v}{\partial \theta^2} \in L^p_{\gamma_0}(S_{\omega,\rho}) \quad \text{and} \quad \frac{\partial^2 v}{\partial r \partial \theta} \in L^p_{\gamma_1}(S_{\omega,\rho}).$$
Now, we explicit the fact that  
$$\frac{\partial^3 V_1}{\partial t \partial\theta^2}, \frac{\partial^3 V_1}{\partial t^2\partial\theta}, \frac{\partial^4 V_1}{\partial t^2\partial\theta^2} \in L^p \left((0,+\infty)\times (0,\omega)\right)$$ 
We have 
$$\frac{\partial^3 V_1}{\partial t^2 \partial \theta}(t,\theta) = \frac{(\nu+1)^2}{\left(\rho e^{-t}\right)^{\nu + 1}}  \frac{\partial v}{\partial \theta}(\rho e^{-t},\theta) - \frac{2\nu + 1}{\left(\rho e^{- t}\right)^{\nu}}  \frac{\partial^2 v}{\partial r \partial\theta}(\rho e^{-t},\theta) + \frac{1}{\left(\rho e^{-t}\right)^{\nu-1}}  \frac{\partial^3 v}{\partial r^2 \partial \theta}(\rho e^{-t},\theta),$$
hence
$$\frac{1}{\left(\rho e^{-t}\right)^{\nu-1}}  \frac{\partial^3 v}{\partial r^2 \partial \theta}(\rho e^{-t},\theta) = \frac{\partial^3 V_1}{\partial t^2 \partial \theta}(t,\theta) - \frac{(\nu+1)^2}{\left(\rho e^{-t}\right)^{\nu + 1}}  \frac{\partial v}{\partial \theta}(\rho e^{-t},\theta) + \frac{2\nu + 1}{\left(\rho e^{- t}\right)^{\nu}}  \frac{\partial^2 v}{\partial r \partial\theta}(\rho e^{-t},\theta).$$
Thus, we obtain
$$\frac{\partial^3 v}{\partial r^2 \partial \theta} \in L^p_{\gamma_2}(S_{\omega,\rho}),$$
and in the same way, we have also
$$\frac{\partial^3 v}{\partial r \partial\theta^2} \in L^p_{\gamma_1}(S_{\omega,\rho}) \quad \text{and} \quad \frac{\partial^4 v}{\partial r^2 \partial \theta^2} \in L^p_{\gamma_2}(S_{\omega,\rho}).$$
Now, we explicit the fact that
$$V_1 \in L^p\left((0,+\infty);W^{4,p}(0,\omega)\right),$$
that is, for all $i = 1,2,3,4$
$$\int_0^{+\infty} \left\|\frac{\partial^i V_1}{\partial \theta^i}(t,\theta)\right\|^p_{L^p(0,\omega)} dt < +\infty.$$
Then, we have
$$\begin{array}{lll}
\dis\int_0^{+\infty} \int_0^\omega \left|\frac{\partial^i V_1}{\partial \theta^i}(t,\theta)\right|^p d\theta~dt & = &\dis\int_0^{+\infty} \int_0^\omega \left|\frac{1}{\left(\rho e^{-t}\right)^{\nu+1}}\frac{\partial^i v}{\partial\theta^i}(\rho e^{-t},\theta)\right|^p d\theta~dt \\ \ecart

& = & \dis \int_0^\rho \int_0^\omega r^{-4p + 1}\left|\frac{\partial^i v}{\partial\theta^i}\left(r,\theta\right)\right|^p d\theta~dr,
\end{array}$$
which gives
$$\frac{\partial^i v}{\partial\theta^i} \in L^p_{\gamma_0}(S_{\omega,\rho}), \quad \text{for }i = 1,2,3,4.$$

\subsubsection{Regularity of $V_2$}

In the same way, we explicit the fact that
$$V_2 \in W^{2,p}\left((0,+\infty)\times (0,\omega)\right),$$
where 
$$V_2 (t,\theta) = \left(r \frac{\partial}{\partial r}\right)^2 \left(\frac{v}{r}\right)(\rho e^{-t},\theta) = r\frac{\partial^2 v}{\partial r^2}(\rho e^{-t},\theta) - \frac{\partial v}{\partial r}(\rho e^{-t},\theta) + \frac{v}{r}(\rho e^{-t},\theta).$$
It is clear, from subsection \ref{SubSect First situation}, that
$$r\frac{\partial^2 v}{\partial r^2},~\frac{\partial v}{\partial r},~\frac{v}{r} \in L^p_{\gamma_1}(S_{\omega,\rho}).$$
Moreover, we have
$$\frac{\partial V_2}{\partial t}(t,\theta) = -r\rho e^{-t} \frac{\partial^3 v}{\partial r^3} (\rho e^{-t}, \theta) + \rho e^{-t} \frac{\partial^2 v}{\partial r^2} (\rho e^{-t}, \theta) - \frac{\rho e^{-t}}{r} \frac{\partial v}{\partial r} (\rho e^{-t}, \theta),$$
hence
$$r\rho e^{-t} \frac{\partial^3 v}{\partial r^3} (\rho e^{-t}, \theta) = \rho e^{-t} \frac{\partial^2 v}{\partial r^2} (\rho e^{-t}, \theta) - \frac{\rho e^{-t}}{r} \frac{\partial v}{\partial r} (\rho e^{-t}, \theta) - \frac{\partial V_2}{\partial t}(t,\theta).$$
Then
$$\begin{array}{lll}
\dis\int_0^{+\infty} \int_0^\omega \left| \rho e^{-t} \frac{\partial^2 v}{\partial r^2} (\rho e^{-t}, \theta) \right|^p d\theta~dt & = &\dis\int_0^{\rho} \int_0^\omega \frac{r^{3p - 1}}{r}r^{-3p+1}\left|r \frac{\partial^2 v}{\partial r^2}(r,\theta)\right|^p d\theta~dr \\ \ecart

& = & \dis \int_0^\rho \int_0^\omega r^{3p - 2} r^{-3p+1}\left|r \frac{\partial^2 v}{\partial r^2}(r,\theta)\right|^p d\theta~dr \\ \ecart

& \leqslant & \dis \rho^{3p - 2} \int_0^\rho \int_0^\omega  r^{-3p+1}\left|r \frac{\partial^2 v}{\partial r^2}(r,\theta)\right|^p d\theta~dr < + \infty,
\end{array}$$
and 
$$\begin{array}{lll}
\dis\int_0^{+\infty} \int_0^\omega \left|\frac{\rho e^{-t}}{r} \frac{\partial v}{\partial r} (\rho e^{-t}, \theta) \right|^p d\theta~dt & = &\dis\int_0^{\rho} \int_0^\omega \frac{r^{3p - 1}}{r}r^{-3p+1}\left|\frac{\partial v}{\partial r}(r,\theta)\right|^p d\theta~dr \\ \ecart

& \leqslant & \dis \rho^{3p - 2} \int_0^\rho \int_0^\omega  r^{-3p+1}\left|\frac{\partial v}{\partial r}(r,\theta)\right|^p d\theta~dr < + \infty.
\end{array}$$
It follows that
$$\int_0^{+\infty} \int_0^\omega \left|r\rho e^{-t} \frac{\partial^3 v}{\partial r^3}(\rho e^{-t},\theta)\right|^p d\theta~dt = \int_0^{\rho} \int_0^\omega r^{2p - 1}\left|\frac{\partial^3 v}{\partial r^3}(r,\theta)\right|^p d\theta~dr < + \infty,$$
which means that 
\begin{equation*}
\frac{\partial^3 v}{\partial r^3} \in L^p_{\gamma_3}(S_{\omega,\rho}),
\end{equation*}
where $\gamma_3 = 2 - \frac{1}{p}$.

In the same way, we have
\begin{equation*}
\frac{\partial^4 v}{\partial r \partial \theta^3} \in L^p_{\gamma_1}(S_{\omega,\rho}) \quad \text{and} \quad \frac{\partial^4 v}{\partial r^3 \partial \theta} \in L^p_{\gamma_3}(S_{\omega,\rho}).
\end{equation*}
Furthermore, we have
$$\begin{array}{lll}
\dis\frac{\partial^2 V_2}{\partial t^2}(t,\theta) & = & \dis r\rho e^{-t} \frac{\partial^3 v}{\partial r^3} (\rho e^{-t}, \theta) + r\left(\rho e^{-t}\right)^2 \frac{\partial^4 v}{\partial r^4} (\rho e^{-t}, \theta) - \rho e^{-t} \frac{\partial^2 v}{\partial r^2} (\rho e^{-t}, \theta) \\ \ecart
&& \dis - \left(\rho e^{-t}\right)^2 \frac{\partial^3 v}{\partial r^3} (\rho e^{-t}, \theta) + \frac{\rho e^{-t}}{r}\frac{\partial v}{\partial r} (\rho e^{-t}, \theta) + \frac{\left(\rho e^{-t}\right)^2}{r}\frac{\partial^2 v}{\partial r^2} (\rho e^{-t}, \theta),
\end{array}$$
hence
$$\int_0^{+\infty} \int_0^\omega \left| r\left(\rho e^{-t}\right)^2 \frac{\partial^4 v}{\partial r^4} (\rho e^{-t}, \theta) \right|^p d\theta~dt = \int_0^{\rho} \int_0^\omega r^{3p - 1}\left|\frac{\partial^4 v}{\partial r^4}(r,\theta)\right|^p d\theta~dr < + \infty.$$
Then
$$\frac{\partial^4 v}{\partial r^4} \in L^p_{\gamma_4}(S_{\omega,\rho}),$$
where $\gamma_4 = 3 - \frac{1}{p}$.

\end{document}